\documentclass[12pt]{amsproc}
\usepackage[cp1251]{inputenc}
\usepackage[T2A]{fontenc} 
\usepackage[russian]{babel} 
\usepackage{amsbsy,amsmath,amssymb,amsthm,graphicx,color} 
\usepackage{amscd}
\def\le{\leqslant}

\newtheorem{prop}{Предложение}

\theoremstyle{definition}

\theoremstyle{remark}

\begin {document}
\unitlength=1mm
\title[Двойственность Саито]
{Двойственность Саито и классическая теория арифметических функций}
\author{Г. Г. Ильюта}
\email{ilyuta@mccme.ru}
\address{}
\thanks{Работа поддержана грантами РФФИ-16-01-00409 и НШ-5138.2014.1.}

\bigskip

\begin{abstract}
Мы изучаем двойственность Саито и преобразование Фурье-Рамануджана для степенных сумм и кратностей корней монодромии.

We study Saito duality and Fourier-Ramanujan transform for power sums and multiplicities of monodromy roots.
\end{abstract}

\maketitle
\tableofcontents

\bigskip

\section{Введение}

\bigskip

  Дзета-функция и характеристический многочлен монодромии особенности представимы в виде 
$$
\zeta_e(q):=\prod_{d|n}(q^d-1)^{e(d)} \eqno (1)
$$
для некоторой целочисленной функции $e$, определённой на положительных делителях $d$ фиксированного во всей статье натурального числа $n$ \cite{1}. Двойственная по Саито функция имеет вид \cite{17}
$$
\hat\zeta_e(q):=\prod_{d|n}(q^d-1)^{-e\left(\frac{n}{d}\right)}. 
$$
Нам будет более удобно рассматривать функцию
$$
\zeta^*_e(q):=\hat\zeta^{-1}_e(q)=\prod_{d|n}(q^d-1)^{e\left(\frac{n}{d}\right)}, 
$$
которую назовём преобразованием Саито функции $\zeta_e(q)$. Мы изучаем действие преобразования Саито на (обобщённые) степенные суммы и кратности корней функции $\zeta_e(q)$, используя тот факт, что эти величины связаны с функцией $e$ преобразованием Мёбиуса. Наш подход основан на следующем замечании. Используя формулы для циклотомических многочленов $\Phi_n(q)$ ($\mu(k)$ -- функция Мёбиуса)
$$
q^n-1=\prod_{d|n}\Phi_d(q),
$$
$$
\Phi_n(q)=\prod_{d|n}(q^d-1)^{\mu\left(\frac{n}{d}\right)},\eqno (2)
$$
рациональную функцию вида (1) можно представить как произведение многочленов $\Phi_d(q)$
$$
\prod_{d|n}(q^d-1)^{e(d)}=\prod_{d|n}\Phi_d^{m(\frac{n}{d})}(q). \eqno (3)
$$
Применяя логарифмическую производную к обеим частям этого равенства, получим утверждение о равномерных (even) арифметических функциях \cite{4} -- определённых на множестве целых чисел функциях $a(k)$, зависящих только от наибольшего общего делителя $(k,n)$. А именно, в правой части формулы (3) равномерной является функция $m(k)$ -- кратность корня $e^{\frac{2\pi ik}{n}}$ функции $\zeta_e(q)$, а в левой части после применения логарифмической производной появляется равномерная функция $p(k)$ -- сумма $k$-х степеней корней функции $\zeta_e(q)$. Приведём некоторые утверждения о равномерных функциях (см. \cite{4}, \cite{5}, \cite{12}, \cite{20}), которые будут использоваться ниже. Любая равномерная функция $a(k)$ имеет представление вида (формулы (9) и (10) для $m(k)$ и $p(k)$)
$$
a(k)=\sum_{d|(k,n)}e(d),  k\in\mathbb Z,\eqno (4)
$$
для некоторой функции $e(d)$, $d|n$, а также разложение Фурье-Рамануджана и его обращение
$$
a(k)=\sum_{d|n}r\left(\frac{n}{d}\right)c_d(k),  k\in\mathbb Z,\eqno (5)
$$
$$
r(k)=\frac{1}{n}\sum_{d|n}a\left(\frac{n}{d}\right)c_d(k),  k\in\mathbb Z,\eqno (6)
$$
другими словами, суммы Рамануджана 
$$
c_m(l)=\sum\limits_{\substack{k=1\\(k,m)=1}}^me^{\frac{2\pi ikl}{m}}, l\in\mathbb Z,
$$
для $m|n$ образуют базис в пространстве равномерных функций. Заметим, что равномерная функция $a(k)$ из формулы (4) представляет собой кратность корня $1$ частичной дзета-функции
$$
\zeta_e^{(k)}(q):=\sum_{d|(k,n)}(q^d-1)^{e(d)},  k\in\mathbb Z.
$$
Из равенств (5), (6) и формулы Рамануджана
$$
\sum_{k=1}^\infty \frac{c_n(k)}{k^s}=\zeta(s)\phi_{1-s}(n),
$$
где
$$
\phi_s(n):=\sum_{d|n}\mu\left(\frac{n}{d}\right)d^s,
$$
вытекают равенства
$$
\sum_{k=1}^\infty \frac{a(k)}{k^s}=\zeta(s)\sum_{d|n}r\left(\frac{n}{d}\right)\phi_{1-s}(d),\eqno (7)
$$
$$
\sum_{k=1}^\infty \frac{r(k)}{k^s}=\frac{\zeta(s)}{n}\sum_{d|n}a\left(\frac{n}{d}\right)\phi_{1-s}(d).\eqno (8)
$$

   В $\S 2$ получены равенства для различных производящих функций последовательностей $m(k)$, $p(k)$ и $e(d)$. В частности, мы свяжем многочлены Апостола-Бернулли и Апостола-Эйлера с производящими функциями, которые являются $q$-степенными суммами показателей Кокстера в случае простых особенностей. В \cite{3} степенные суммы показателей Кокстера представлены как специализации многочленов Тодда. В \cite{16} К.~ Саито упоминает свой препринт 1985 года, в котором степенные суммы показателей регулярных весовых систем представлены как многочлены от весов. Степенным суммам обобщённых показателей Кокстера (спектральных чисел особенности) посвящены статьи \cite{2}, \cite{8}. 

  При доказательстве формул из $\S 2$ используются равенства (7), (8) и следующее замечание: если две линейные комбинации линейно независимых функций равны, то равенство сохранится при замене этих функций на любые другие. Возможен другой единый подход к доказательству формул из $\S 3$: при умножении (делении) функций вида (1) соответствующие равенства складываются (вычитаются) и поэтому достаточно привести доказательство для $\zeta_e(q)=1-q^d$, $d|n$. Например, число $e^{\frac{2\pi ik}{n}}$ является корнем многочлена $q^{n/d}-1$, $d|n$, если и только если $d|(k,n)$. Поэтому
$$
m(k)=\sum_{d|(k,n)}e\left(\frac{n}{d}\right),  k\in\mathbb Z.\eqno (9)
$$
Другой пример: по определению
$$
p(k)=\sum_{d|n}e(d)\sum_{l=0}^{d-1}e^{\frac{2\pi ilk}{d}}, k\in\mathbb Z,
$$
и поэтому из соотношений ортогональности для характеров циклической группы $C_d$, $d|n$, вытекает равенство
$$
p(k)=\sum_{d|(k,n)}de(d),  k\in\mathbb Z.\eqno (10)
$$
Для функции $\zeta^*_e(q)$
$$
m^*(k):=\sum_{d|(k,n)}e(d),  k\in\mathbb Z,
$$
$$
p^*(k):=\sum_{d|(k,n)}de\left(\frac{n}{d}\right),  k\in\mathbb Z.
$$ 

  В $\S 3$ получены формулы перехода, связывающие обобщения функций $m(k)$, $p(k)$ и их преобразований Саито $m^*(k)$, $p^*(k)$. Обобщения определяются с помощью производящих функций для преобразования Мёбиуса -- функция $(1-q)^{-1}$ и дзета-функция Римана $\zeta(s)$ в этих производящих функциях заменяются, соответственно, формальным степенным рядом и формальным рядом Дирихле. Такие обобщения позволяют учитывать в степенных суммах и кратностях только корни с предписанными свойствами. Отметим ещё одно направление для обобщений: в преобразовании Мёбиуса можно рассматривать суммирование по подмножеству делителей, например, регулярные свёртки Наркиевича \cite{12}.

   В $\S 4$ содержится следствие формулы (3) для эта-произведения Саито. Для квазиоднородных особенностей в $\S 5$ приведены формулы для производящих функций последовательностей $m(k)$ и $p(k)$ (используется формула для функции $\zeta_e(q)$ из \cite{7} и формула для многочлена $\sum_{k=0}^{n-1}m(k)q^k$ из \cite{22}).

  Корнями характеристического многочлена монодромии являются числа $e^{\frac{2\pi im_j}{n}}$, где $m_j\in\mathbb Z$ -- подходящим образом нормализованные спектральные числа особенности \cite{15}. Отсюда следует, что
$$
\sum_{k=0}^{n-1}m(k)q^k=\sum_j q^{m_j} \mod(q^n-1),
$$
в частности, для простых особенностей 
$$
\sum_{k=0}^{n-1}m(k)q^k=\sum_j q^{m_j} 
$$
и числа $m_j$ совпадают с показателями Кокстера соответствующих конечных групп отражений. В $\S 6$ и $\S 7$ для простых, параболических и исключительных унимодальных особенностей приведены таблицы производящих функций последовательностей $m(k)$ и $p(k)$ (при этом используются таблицы функций $e(d)$ из \cite{17}).

  В \cite{9} суммы Рамануджана интерпретированы как элементы таблицы суперхарактеров циклической группы $C_n$ (любая теория суперхарактеров конечной группы определяются специального вида разложением регулярного характера, обычная теория характеров отвечает разложению на неприводимые). Равномерные функции являются функциями суперклассов сопряжённости в этой теории. В контексте формулы А'Кампо для дзета-функции циклическая группа $C_n$ порождается оператором монодромии. В последние годы изучаются и более общие группы симметрий особенностей \cite{8}. Интересно было бы понять роль суперхарактеров этих групп симметрий. 
 
  Мы предполагаем, что областью значений функции $e(d)$ является множество целых чисел $\mathbb Z$, но фактически нам понадобится только возможность рассматривать формальные линейные комбинации значений функции $e(d)$. Например, кратности и степенные суммы корней дзета-функции можно определить формулами (9) и (10), не заботясь о том, что собой представляют переменные $e(d)$. В формуле А'Кампо для дзета-функции монодромии \cite{1} числа $e(d)$ являются эйлеровыми характеристиками некоторых многообразий. В мотивных обобщениях дзета-функции А'Кампо функция $e(d)$ принимает значения в универсальном для эйлеровой характеристики объекте -- кольце Гротендика алгебраических многообразий.

  Записывая степенные суммы корней функции $\zeta_e(q)$ как суммы степеней всех корней степени $n$ из единицы с кратностями, мы получим определение дискретного преобразования Фурье, которым связаны функции $p(k)$ и $m(k)$,
$$
p(l)=\sum_{k=0}^{n-1}m(k)e^{\frac{2\pi ilk}{n}}
$$
$$
=\sum_{d|n}\sum_{(k,n)=d}m(k)e^{\frac{2\pi ilk/d}{n/d}} 
$$
$$
=\sum_{d|n}m(d)\sum_{(k/d,n/d)=1}e^{\frac{2\pi ilk/d}{n/d}} 
$$
$$
=\sum_{d|n}m(d)c_{n/d}(l)=\sum_{d|n}m(n/d)c_d(l), l\in\mathbb Z.  
$$
Заметим, что эти равенства равносильны равенству (3) -- суммы Рамануджана $c_d(l)$, $l\in\mathbb Z$, являются степенными суммами корней многочлена $\Phi_d(q)$. Выбирая произвольно функцию $e$, мы видим, что полученные ниже формулы справедливы для любой пары равномерных функций, которые связаны дискретным преобразованием Фурье. Любая такая пара включается в семейство пар \cite{5}: если
$$
f_s(k):=\sum_{d|(k,n)}\frac{F(d,n/d)}{d^s}, f'_s(k):=\sum_{d|(k,n)}\frac{F(n/d,d)}{(n/d)^s},
$$
то
$$
f_s(k)=\sum_{d|n}f'_{s+1}(n/d)c_d(k),  k\in\mathbb Z.
$$
Преобразование Саито
$$
m(k)\to m^*(k): e(n/d)\to e(d)
$$
и преобразование Фурье-Рамануджана (ограничение на пространство равномерных функций дискретного преобразования Фурье)
$$
m(k)\to p(k): e(n/d)\to de(d)
$$
включаются в семейство преобразований
$$
e(n/d)\to\frac{e(d)}{d^s}.
$$
Полученные ниже формулы легко обобщаются для этих семейств как формулы перехода, связывающие различные их элементы (на уровне производящих функций Дирихле такие переходы осуществляются с помощью подходящих сдвигов аргумента в рядах Дирихле, ниже мы используем сдвиги на $1$). Например, можно использовать собранные в \cite{10} разложения в ряды Дирихле для функций $G_1(s-k)/G_2(s)$ в Примерах 10-21, где используются разложения для функций $G_1(s-1)/G_2(s)$. Мы не будем этим заниматься подробно, но приведём следующий факт, указывающий на возможную связь с теоремой Себастьяни-Тома: характеристический многочлен монодромии прямой суммы особенностей равен тензорному произведению характеристических многочленов слагаемых. В формуле А'Кампо для дзета-функции монодромии \cite{1} многочлен $q^d-1$ появляется как характеристический многочлен циклической $d$-перестановки.

\begin{prop}\label{prop1} Характеристический многочлен $k$-й тензорной степени циклической $d$-перестановки равен
$$
(q^d-1)^{d^{k-1}}.
$$
\end{prop}

  Доказательство. Тензорному произведению линейных операторов соответствует тензорное произведение их характеристических многочленов $f\otimes g(q)$ -- многочлен, корнями которого являются произведения корней многочленов $f(q)$ и $g(q)$.  Хорошо известно, что тензорное произведение связано с результантом формулой
$$
f\otimes g(q)=Res_t(t^mf\left(\frac{q}{t}\right),g(t)), m=\deg f.
$$
Для $k$ тензорных сомножителей по индукции имеем
$$
(q^d-1)\otimes\dots\otimes (q^d-1)=Res_t((q^d-t^d)^{d^{k-2}},q^d-1)
$$
$$
=\prod_{l=0}^{d-1}(q^d-(e^{\frac{2\pi il}{d}})^d)^{d^{k-2}}=(q^d-1)^{d^{k-1}}.
$$

\bigskip

\section{Производящие функции для равномерных функций}

\bigskip

  Известны следующие способы представления преобразования Мёбиуса $a(k)=\sum_{d|k}e(d)$, $k\in\mathbb N$, в виде соотношений между производящими функциями: представление ряда Ламберта в виде степенного ряда
$$
\sum_{k=1}^\infty a(k)q^k=\sum_{m=1}^\infty \frac{e(m)q^m}{1-q^m} \eqno (11)
$$
и соотношение между рядами Дирихле 
$$
\sum_{k=1}^\infty \frac{a(k)}{k^s}=\zeta(s)\sum_{m=1}^\infty \frac{e(m)}{m^s},\eqno (12)
$$
показывающее, что преобразование Мёбиуса последовательности отвечает умножению на дзета-функцию Римана соответствующей производящей функции Дирихле. Перепишем формулу (11) в более близком к формуле (12) виде
$$
\sum_{k=1}^\infty a(k)q^k=\frac{1}{1-q}\sum_{m=1}^\infty \frac{e(m)q^m}{[m]_q},
$$
где $[m]_q=1+q+\dots+q^{m-1}$. Любая равномерная функция (7) является $n$-периодической и поэтому для неё формулы принимают вид
$$
\sum_{k=1}^\infty a(k)q^k=\frac{\sum_{k=1}^n a(k)q^k}{1-q^n}=\sum_{d|n}\frac{e(d)q^d}{1-q^d}=\frac{1}{1-q}\sum_{d|n}\frac{e(d)q^d}{[d]_q},\eqno (13)
$$
$$
\sum_{k=1}^\infty \frac{a(k)}{k^s}=\zeta(s)\sum_{d|n}\frac{e(d)}{d^s}.\eqno (14)
$$

  Многочлены Апостола-Бернулли $B_k(n,q)$ и Апостола-Эйлера $E_k(n,q)$ определяются производящими функциями \cite{21}
$$
\frac{te^{tx}}{qe^t-1}=\sum_{n=0}^\infty B_n(x,q)\frac{t^n}{n!},
$$
$$
\frac{2e^{tx}}{qe^t+1}=\sum_{n=0}^\infty E_n(x,q)\frac{t^n}{n!}.
$$

\begin{prop}\label{prop2} Если
$$
a(k)=\sum_{d|(k,n)}e(d),  k\in\mathbb Z,
$$
то
$$
\frac{\sum_{k=0}^{n-1}a(k)q^k}{1-q^n}=\sum_{d|n}\frac{e(d)}{1-q^d}=\frac{1}{1-q}\sum_{d|n}\frac{e(d)}{[d]_q},\eqno (15)
$$
$$
\sum_{k=0}^{n-1}a(k)(bk+c)^rq^k=\sum_{d|n}\frac{e(d)(bd)^r}{r+1}\left(q^nB_{r+1}\left(\frac{c+bn}{bd},q^d\right)-B_{r+1}\left(\frac{c}{bd},q^d\right)\right),     
$$
$$
\sum_{k=0}^{n-1}a(k)(bk+c)^r(-q)^k
$$
$$
=\sum_{2\nmid d|n}e(d)\frac{(bd)^r}{2}\left((-1)^{\frac{n}{d}-1}q^nE_r\left(\frac{c+bn}{bd},q^d\right)-E_r\left(\frac{c}{bd},q^d\right)\right)
$$
$$
+\sum_{2|d|n}\frac{e(d)(bd)^r}{r+1}\left(q^nB_{r+1}\left(\frac{c+bn}{bd},q^d\right)-B_{r+1}\left(\frac{c}{bd},q^d\right)\right).     
$$
\end{prop}

  Доказательство. Используя равенство $a(0)=a(n)=\sum_{d|n}e(d)$ и формулу (13), получим
$$
\frac{\sum_{k=0}^{n-1}a(k)q^k}{1-q^n}-\sum_{d|n}e(d)=\frac{-a(0)+a(n)q^n+\sum_{k=0}^{n-1}a(k)q^k}{1-q^n}
$$
$$
=\frac{\sum_{k=1}^na(k)q^k}{1-q^n}=\sum_{d|n}\frac{e(d)q^d}{1-q^d}
$$
$$
=\sum_{d|n}\frac{e(d)}{1-q^d}-\sum_{d|n}e(d).
$$

  Перепишем равенство (15) в следующем виде
$$
\sum_{k=0}^{n-1}a(k)q^k=\sum_{d|n}e(d)(1+q^d+\dots+q^{n-d}),
$$
Заменяя $q^k\to (bk+c)^rq^k$ или $q^k\to (bk+c)^r(-q)^k$ для всех $k$, приходим к равенствам
$$
\sum_{k=0}^{n-1}a(k)(bk+c)^rq^k=\sum_{d|n}e(d)\sum_{i=0}^{\frac{n}{d}-1}(bdi+c)^rq^{di},
$$
$$
\sum_{k=0}^{n-1}a(k)(bk+c)^r(-q)^k=\sum_{2\nmid d|n}e(d)\sum_{i=0}^{\frac{n}{d}-1}(-1)^{i}(bdi+c)^rq^{di}
$$
$$
+\sum_{2|d|n}e(d)\sum_{i=0}^{\frac{n}{d}-1}(bdi+c)^rq^{di}.
$$
Остаётся использовать формулы \cite{21}
$$
\sum_{i=0}^n (bi+c)^rq^i=\frac{b^r}{r+1}\left(q^{n+1}B_{r+1}\left(\frac{c}{b}+n+1,q\right)-B_{r+1}\left(\frac{c}{b},q\right)\right),
$$
$$
\sum_{i=0}^n (-1)^i(bi+c)^rq^i=\frac{b^r}{2}\left((-1)^nq^{n+1}E_r\left(\frac{c}{b}+n+1,q\right)-E_r\left(\frac{c}{b},q\right)\right).
$$

  Из формул (7), (8) и формулы (14) для $a(k)=m(k)$ и $a(k)=p(k)$ вытекает

\begin{prop}\label{prop3}
$$
\sum_{d|n}m\left(\frac{n}{d}\right)\phi_{1-s}(d)=\sum_{d|n}\frac{de(d)}{d^s},
$$
$$
\frac{1}{n}\sum_{d|n}p\left(\frac{n}{d}\right)\phi_{1-s}(d)=\sum_{d|n}\frac{e(n/d)}{d^s} \eqno (16)
$$
или, после замены $s\to s+1$,
$$
\sum_{d|n}m\left(\frac{n}{d}\right)\phi_{-s}(d)=\sum_{d|n}\frac{e(d)}{d^s},\eqno (17)
$$
$$
\sum_{d|n}p\left(\frac{n}{d}\right)\phi_{-s}(d)=\sum_{d|n}\frac{n/de(n/d)}{d^s}.
$$
\end{prop}
 
  Предложения 4-6 являются следствиями Предложения 3 (последовательность $\phi_{-s}(d)$, $d|n$, является обратным преобразованием Мёбиуса последовательности $d^{-s}$, $d|n$).

\begin{prop}\label{prop4} Если последовательность $z_d$, $d|n$, является обратным преобразованием Мёбиуса последовательности $x_d$, $d|n$,
$$
z_d=\sum_{d'|d}\mu\left(\frac{d}{d'}\right)x_{d'},
$$
то
$$
\sum_{d|n}m\left(\frac{n}{d}\right)z_d=\sum_{d|n}e(d)x_d,
$$
$$
\sum_{d|n}p\left(\frac{n}{d}\right)z_d=\sum_{d|n}n/de(n/d)x_d.
$$
\end{prop}

  Обратным преобразованием Мёбиуса последовательности $q^d$, $d|n$, является последовательность $dM(q,d)$, $d|n$, где $M(q,d)$ -- многочлен ожерелий \cite{13},
$$
M(q,d)=\frac{1}{d}\sum_{d'|d}\mu\left(\frac{d}{d'}\right)q^{d'}.
$$

\begin{prop}\label{prop5}
$$
\sum_{d|n}m\left(\frac{n}{d}\right)dM(q,d)=\sum_{d|n}e(d)q^d,
$$ 
$$
\sum_{d|n}p\left(\frac{n}{d}\right)dM(q,d)=\sum_{d|n}n/de(n/d)q^d.
$$
\end{prop}

  Применяя логарифмическую производную к обеим частям равенства (2), получим
обратное преобразование Мёбиуса
$$
\frac{q\Phi'_n(q)}{\Phi_n(q)}=\sum_{d|n}\mu\left(\frac{n}{d}\right)\frac{dq^d}{q^d-1}.
$$

\begin{prop}\label{prop6}
$$
\sum_{d|n}m\left(\frac{n}{d}\right)\frac{q\Phi'_d(q)}{\Phi_d(q)}=\sum_{d|n}\frac{de(d)q^d}{q^d-1},
$$ 
$$
\sum_{d|n}p\left(\frac{n}{d}\right)\frac{q\Phi'_d(q)}{\Phi_d(q)}=\sum_{d|n}\frac{ne(n/d)q^d}{q^d-1}.
$$
\end{prop}

Согласно \cite{14}
$$
\frac{q\Phi'_n(q)}{\Phi_n(q)}=\frac{\sum_{k=1}^nc_n(k)q^k}{q^n-1}.\eqno (18)
$$
Поэтому из Предложения 6 следует
\begin{prop}\label{prop7}
$$
\sum_{d|n}m\left(\frac{n}{d}\right)\frac{\sum_{k=1}^dc_d(k)q^k}{q^d-1}=\sum_{d|n}\frac{de(d)q^d}{q^d-1},
$$ 
$$
\sum_{d|n}p\left(\frac{n}{d}\right)\frac{\sum_{k=1}^dc_d(k)q^k}{q^d-1}=\sum_{d|n}\frac{ne(n/d)q^d}{q^d-1}.
$$
\end{prop}

\bigskip

\section{Обобщённые степенные суммы и кратности корней}

\bigskip

Обобщая формулу (14), для формального ряда Дирихле $G(s)=\sum_{k=1}^\infty g(k)k^{-s}$ определим $m_{G}(k)$, $p_{G}(k)$, $m^*_{G}(k)$, $p^*_{G}(k)$ равенствами
$$
\sum_{k=1}^\infty\frac{m_{G}(k)}{k^s}=G(s)\sum_{d|n}\frac{e\left(\frac{n}{d}\right)}{d^s},
$$
$$
\sum_{k=1}^\infty\frac{p_{G}(k)}{k^s}=G(s)\sum_{d|n}\frac{de(d)}{d^s}.
$$
$$
\sum_{k=1}^\infty\frac{m^*_{G}(k)}{k^s}=G(s)\sum_{d|n}\frac{e(d)}{d^s},
$$
$$
\sum_{k=1}^\infty\frac{p^*_{G}(k)}{k^s}=G(s)\sum_{d|n}\frac{de\left(\frac{n}{d}\right)}{d^s}.
$$
Аналогично можно обобщить формулу (15)
$$
\sum_{k=1}^\infty m_{[G]}(k)q^k=\sum_{k=1}^\infty g(k)q^k\sum_{d|n}\frac{e\left(\frac{n}{d}\right)}{[d]_q},
$$
$$
\sum_{k=1}^\infty p_{[G]}(k)q^k=\sum_{k=1}^\infty g(k)q^k\sum_{d|n}\frac{de(d)}{[d]_q}.
$$

\begin{prop}\label{prop8} 
$$
G(s)\sum_{d|n}m\left(\frac{n}{d}\right)\phi_{-s}(d)=\sum_{k=1}^\infty\frac{m^*_{G}(k)}{k^s}.
$$
$$
\frac{1}{n}G(s)\sum_{d|n}p\left(\frac{n}{d}\right)\phi_{2-s}(d)=\sum_{k=1}^\infty\frac{p^*_{G}(k)}{k^s}.
$$
\end{prop}

  Доказательство. Умножаем формулу (17) на $G(s)$ и используем определение $m^*_{G}(k)$. Формулу (16) после замены $s\to s-1$ умножаем на $G(s)$ и используем определение $p^*_{G}(k)$.

\begin{prop}\label{prop9} 
$$
G_2(s)\sum_{k=1}^\infty\frac{km_{G_1}(k)}{k^s}=G_1(s-1)\sum_{k=1}^\infty\frac{p^*_{G_2}(k)}{k^s},
$$
\end{prop} 

  Доказательство. 
$$
G_2(s)\sum_{k=1}^\infty\frac{km_{G_1}(k)}{k^s}=G_2(s)\sum_{k=1}^\infty\frac{m_{G_1}(k)}{k^{s-1}}
$$
$$
=G_2(s)G_1(s-1)\sum_{d|n}\frac{e\left(\frac{n}{d}\right)}{d^{s-1}}=G_1(s-1)G_2(s)\sum_{d|n}\frac{de\left(\frac{n}{d}\right)}{d^s}
$$
$$
=G_1(s-1)\sum_{k=1}^\infty\frac{p^*_{G_2}(k)}{k^s}.
$$

  Использованные в Примерах 1-12 формулы представлены в каталоге рядов Дирихле \cite{10}. С помощью многих других рядов Дирихле из этого каталога и из \cite{12} можно получить аналогичные следствия Предложения 9.

  Пример 1.
$$
G_1(s)=G_2(s)=\zeta(s).
$$
Используя известные формулы \cite{10}
$$
\frac{1}{\zeta(s)}=\sum_{k=1}^\infty\frac{\mu(k)}{k^s},
$$
$$
\frac{\zeta(s-1)}{\zeta(s)}=\sum_{k=1}^\infty\frac{\phi(k)}{k^s},
$$
$$
\frac{\zeta(s)}{\zeta(s-1)}=\sum_{k=1}^\infty\frac{\phi^{(-1)}(k)}{k^s},
$$
где $\phi(k)$ -- функция Эйлера и
$$
\phi^{(-1)}(k)=\sum_{d|k}d\mu(d)=\sum_{prime p|k}(1-p),
$$
получим
$$
\sum_{k=1}^\infty\frac{km(k)}{k^s}=\left(\sum_{k=1}^\infty\frac{\phi(k)}{k^s}\right)\left(\sum_{k=1}^\infty\frac{p^*(k)}{k^s}\right),
$$
$$
m(k)=\frac{1}{k}\sum_{d|k}\phi\left(\frac{k}{d}\right)p^*(d)
$$
$$
=\frac{1}{k}\sum_{d|k}\sum_{(j/d,k/d)=1}p^*(d)=\frac{1}{k}\sum_{d|k}\sum_{(j,k)=d}p^*(d)
$$
$$
=\frac{1}{k}\sum_{j=1}^kp^*((j,k)), k\in\mathbb N,
$$
$$
\sum_{k=1}^\infty\frac{p^*(k)}{k^s}=\left(\sum_{k=1}^\infty\frac{\phi^{(-1)}(k)}{k^s}\right)\left(\sum_{k=1}^\infty\frac{km(k)}{k^s}\right),
$$
$$
p^*(k)=\sum_{d|k}\phi^{(-1)}\left(\frac{k}{d}\right)dm(d), k\in\mathbb N.
$$

  Пример 2. Для фиксированного $r\in\mathbb N$
$$
G_1(s)=\sum_{d|r}\frac{1}{d^s}, G_2(s)=\zeta(s).
$$
Заметим, что для $r=0$ получаем Пример 1. Согласно \cite{12}, Th. 5.5.,
$$
\frac{G_1(s-1)}{\zeta(s)}=\sum_{k=1}^\infty\frac{c_k(r)}{k^s}
$$
и поэтому
$$
\sum_{k=1}^\infty\frac{km_{G_1}(k)}{k^s}=\left(\sum_{k=1}^\infty\frac{c_k(r)}{k^s}\right)\left(\sum_{k=1}^\infty\frac{p^*(k)}{k^s}\right),
$$
$$
m_{G_1}(k)=\frac{1}{k}\sum_{d|k}c_{k/d}(r)p^*(d), k\in\mathbb N.
$$

  Пример 3. Для фиксированного $r\in\mathbb N$
$$
G_1(s)=\zeta(s), G_2(s)=\frac{1}{\zeta(rs)}.
$$
Согласно \cite{12}, p. 229,
$$
\zeta(s-1)\zeta(rs)=\sum_{k=1}^\infty\frac{\rho_r(k)}{k^s}
$$
где $\rho_r(k)$ равна сумме делителей $d$ числа $k$, для которых число $k/d$ является $r$-й степенью. Поэтому
$$
\sum_{k=1}^\infty\frac{km(k)}{k^s}=\left(\sum_{k=1}^\infty\frac{\rho_r(k)}{k^s}\right)\left(\sum_{k=1}^\infty\frac{p^*_{G_2}(k)}{k^s}\right),
$$
$$
m(k)=\frac{1}{k}\sum_{d|k}\rho_r\left(\frac{k}{d}\right)p^*_{G_2}(d), k\in\mathbb N.
$$

  Пример 4. Для фиксированного $r\in\mathbb N$
$$
G_1(s)=\zeta(s), G_2(s)=\zeta(rs).
$$
Согласно \cite{18}
$$
\frac{\zeta(s-1)}{\zeta(rs)}=\sum_{k=1}^\infty\frac{\Phi_r(k)}{k^s},
$$
где $\Phi_r(k)$ -- функция Кли, равная количеству целых чисел $j$, $1\le j\le k$, для которых $(j,k)_r=1$, $(j,k)_r$ -- наибольший общий делитель, являющийся $r$-й степенью. Поэтому
$$
\sum_{k=1}^\infty\frac{km(k)}{k^s}=\left(\sum_{k=1}^\infty\frac{\Phi_r(k)}{k^s}\right)\left(\sum_{k=1}^\infty\frac{p^*_{G_2}(k)}{k^s}\right),
$$
$$
m(k)=\frac{1}{k}\sum_{d|k}\Phi_r\left(\frac{k}{d}\right)p^*_{G_2}(d), k\in\mathbb N.
$$

  Пример 5. 
$$
G_1(s)=\zeta(s), G_2(s)=\frac{\zeta(s)}{\zeta(2s)}.
$$
Согласно \cite{11}, p. 255,
$$
\frac{\zeta(s)}{\zeta(2s)}=\sum_{k=1}^\infty\frac{|\mu(k)|}{k^s}
$$
и согласно \cite{12}, p.229,
$$
\frac{\zeta(s-1)\zeta(2s)}{\zeta(s)}=\sum_{k=1}^\infty\frac{\beta(k)}{k^s},
$$
где $\beta(k)$ равна количеству целых чисел $j$, $1\le j\le k$, для которых $(j,k)$ является квадратом. Поэтому
$$
\sum_{k=1}^\infty\frac{km(k)}{k^s}=\left(\sum_{k=1}^\infty\frac{\beta(k)}{k^s}\right)\left(\sum_{k=1}^\infty\frac{p^*_{G_2}(k)}{k^s}\right),
$$
$$
m(k)=\frac{1}{k}\sum_{d|k}\beta\left(\frac{k}{d}\right)p^*_{G_2}(d), k\in\mathbb N.
$$

  Пример 6. 
$$
G_1(s)=\zeta(s), G_2(s)=\frac{\zeta(2s)}{\zeta(s)}.
$$
Согласно \cite{12}, p.230,
$$
\frac{\zeta(s-1)\zeta(s)}{\zeta(2s)}=\sum_{k=1}^\infty\frac{\psi(k)}{k^s},
$$
где $\psi(k)=\sum_{d|k}|\mu(k/d)|$ -- функция Дедекинда. Поэтому
$$
\sum_{k=1}^\infty\frac{km(k)}{k^s}=\left(\sum_{k=1}^\infty\frac{\psi(k)}{k^s}\right)\left(\sum_{k=1}^\infty\frac{p^*_{G_2}(k)}{k^s}\right),
$$
$$
m(k)=\frac{1}{k}\sum_{d|k}\psi\left(\frac{k}{d}\right)p^*_{G_2}(d), k\in\mathbb N.
$$

  Пример 7. 
$$
G_1(s)=G_2(s)=\frac{\zeta(2s)}{\zeta(s)}.
$$
Согласно \cite{12}, p. 232,
$$
\frac{\zeta(s)\zeta(2(s-1))}{\zeta(s-1)\zeta(2s)}=\sum_{k=1}^\infty\frac{\lambda(k)\phi(k)}{k^s},
$$
где $\lambda(k)=(-1)^{\sum n_i}$ -- функция Лиувилля, $k=\prod p_i^{n_i}$ -- разложение на простые. Поэтому
$$
\sum_{k=1}^\infty\frac{km_{G_1}(k)}{k^s}=\left(\sum_{k=1}^\infty\frac{\lambda(k)\phi(k)}{k^s}\right)\left(\sum_{k=1}^\infty\frac{p^*_{G_2}(k)}{k^s}\right),
$$
$$
m_{G_1}(k)=\frac{1}{k}\sum_{d|k}\lambda\left(\frac{k}{d}\right)\phi\left(\frac{k}{d}\right)p^*_{G_2}(d), k\in\mathbb N.
$$

  Пример 8. Для фиксированного $r\in\mathbb N$
$$
G_1(s)=\frac{\zeta(2rs)}{\zeta(rs)}, G_2(s)=\frac{\zeta(s)}{\zeta(2s)}.
$$
Согласно \cite{12}, p. 232,
$$
\frac{\zeta(2s)\zeta(2r(s-1))}{\zeta(s)\zeta(r(s-1))}=\sum_{k=1}^\infty\frac{\lambda(k)\rho'_r(k)}{k^s},
$$
где $\rho'_r(k)$ равна сумме делителей числа $k$, являющихся $r$-й степенью. Поэтому
$$
\sum_{k=1}^\infty\frac{km_{G_1}(k)}{k^s}=\left(\sum_{k=1}^\infty\frac{\lambda(k)\rho'_r(k)}{k^s}\right)\left(\sum_{k=1}^\infty\frac{p^*_{G_2}(k)}{k^s}\right),
$$
$$
m_{G_1}(k)=\frac{1}{k}\sum_{d|k}\lambda\left(\frac{k}{d}\right)\rho'_r\left(\frac{k}{d}\right)p^*_{G_2}(d), k\in\mathbb N.
$$

  Пример 9. Для фиксированного $r\in\mathbb N$
$$
G_1(s)=\sum_{d|n}\frac{(r,d)^{s+1}\mu\left(\frac{n}{d}\right)}{d^s}, G_2(s)=\frac{1}{\zeta(s)}.
$$
Согласно \cite{6}
$$
\zeta(s)G_1(s-1)=\sum_{k=1}^\infty\frac{c_n(rk)}{k^s}
$$
и поэтому
$$
\sum_{k=1}^\infty\frac{km_{G_1}(k)}{k^s}=\left(\sum_{k=1}^\infty\frac{c_n(rk)}{k^s}\right)\left(\sum_{k=1}^\infty\frac{p^*_{G_2}(k)}{k^s}\right),
$$
$$
m_{G_1}(k)=\frac{1}{k}\sum_{d|k}c_n\left(\frac{rk}{d}\right)p^*_{G_2}(d), k\in\mathbb N.
$$

  Пример 10.
$$
G_1(s)=\sum_{d|n}\frac{\lambda(d)\mu\left(\frac{n}{d}\right)}{d^s}, G_2(s)=\frac{\zeta(s)}{\zeta(2s)}.
$$
Согласно \cite{6}
$$
\frac{\zeta(2s)G_1(s-1)}{\zeta(s)}=\sum_{k=1}^\infty\frac{\lambda(k)c_n(k)}{k^s}
$$
и поэтому
$$
\sum_{k=1}^\infty\frac{km_{G_1}(k)}{k^s}=\left(\sum_{k=1}^\infty\frac{\lambda(k)c_n(k)}{k^s}\right)\left(\sum_{k=1}^\infty\frac{p^*_{G_2}(k)}{k^s}\right),
$$
$$
m_{G_1}(k)=\frac{1}{k}\sum_{d|k}\lambda\left(\frac{k}{d}\right)c_n\left(\frac{k}{d}\right)p^*_{G_2}(d), k\in\mathbb N.
$$

  Пример 11.
$$
G_1(s)=\frac{1}{2^s}-1, G_2(s)=\frac{1}{\zeta(s)}.
$$
Согласно \cite{19}, p. 21,
$$
\zeta(s)G_1(s-1)=\sum_{k=1}^\infty\frac{(-1)^k}{k^s}
$$
и поэтому
$$
\sum_{k=1}^\infty\frac{km_{G_1}(k)}{k^s}=\left(\sum_{k=1}^\infty\frac{(-1)^k}{k^s}\right)\left(\sum_{k=1}^\infty\frac{p^*_{G_2}(k)}{k^s}\right),
$$
$$
m_{G_1}(k)=\frac{1}{k}\sum_{d|k}(-1)^{k/d}p^*_{G_2}(d), k\in\mathbb N.
$$

  Пример 12.
$$
G_1(s)=\zeta(s)(1-\frac{1}{2^s}), G_2(s)=1-\frac{1}{2^s}.
$$
Согласно \cite{19}, p. 6,
$$
\zeta(s-1)\frac{1-\frac{1}{2^{s-1}}}{1-\frac{1}{2^s}}=\sum_{k=1}^\infty\frac{c(k)}{k^s}
$$
где $c(k)$ -- наибольший нечётный делитель. Поэтому
$$
\sum_{k=1}^\infty\frac{km_{G_1}(k)}{k^s}=\left(\sum_{k=1}^\infty\frac{c(k)}{k^s}\right)\left(\sum_{k=1}^\infty\frac{p^*_{G_2}(k)}{k^s}\right),
$$
$$
m_{G_1}(k)=\frac{1}{k}\sum_{d|k}c\left(\frac{k}{d}\right)p^*_{G_2}(d), k\in\mathbb N.
$$

\bigskip

\section{Эта-произведение Саито}

\bigskip

  Эта-произведение Саито $\eta_e(q)$ для функции $e(d)$, $d|n$, определяется формулой \cite{17}
$$
q^{-\frac{\mu_e}{24}}\eta_e(q)=\prod_{k=1}^\infty\zeta_e(q^k):=\bar\eta_e(q),
$$
где $\mu_e=\sum_{d|n}e(d)$.

\begin{prop}\label{prop14}
$$
q\frac{d}{dq}\log\bar\eta_e(q)=\sum_{d|n}de(d)L(q^d)
$$
$$
=\sum_{d|n}m\left(\frac{n}{d}\right)\sum_{k=1}^\infty\frac{kq^k\Phi'_d(q^k)}{\Phi_d(q^k)}
$$
$$
=\sum_{d|n}m\left(\frac{n}{d}\right)\sum_{k=1}^\infty\frac{\sum_{j=1}^dc_d(j)q^{kj}}{q^{kd}-1},
$$
где
$$
L(q):=\sum_{k=1}^\infty\frac{kq^k}{1-q^k}=\sum_{m=1}^\infty\sigma(m)q^m,
$$
$$
\sigma(m):=\sum_{d|m}d.
$$
\end{prop}

  Доказательство. Из формулы (3) следует, что
$$
\bar\eta_e(q)=\prod_{k=1}^\infty\prod_{d|n}(1-q^{kd})^{e(d)}=\prod_{k=1}^\infty\prod_{d|n}\Phi_d^{m(\frac{n}{d})}(q^k),
$$
Применяем оператор $q\frac{d}{dq}\log$ и используем формулу (18).

\bigskip

\section{Квазиоднородные особенности}

\bigskip

Для квазиоднородной особенности функции трёх переменных с регулярной системой весов $(a,b,c;n)$ известна формула для производящей функции спектральных чисел \cite{17}
$$
\sum_j q^{m_j}=q^{-n}\frac{(q^n-q^a)(q^n-q^b)(q^n-q^c)}{(q^a-1)(q^b-1)(q^c-1)}.
$$
Согласно \cite{22} для кратностей $m(k)$ имеем соотношение
$$
\frac{\sum_{k=0}^{n-1}m(k)q^k}{q^n-1}=\frac{n^2}{abc}\frac{1}{q-1}-\frac{1}{q^n-1}      
$$
$$
+\frac{(a,n)}{a}\frac{1}{q^{(a,n)}-1}+\frac{(b,n)}{b}\frac{1}{q^{(b,n)}-1}+\frac{(c,n)}{c}\frac{1}{q^{(c,n)}-1}   
$$
$$
-\frac{n(b,c,n)}{bc}\frac{1}{q^{(b,c,n)}-1}-\frac{n(a,c,n)}{ac}\frac{1}{q^{(a,c,n)}-1}-\frac{n(a,b,n)}{ab}\frac{1}{q^{(a,b,n)}-1}.   
$$
Поэтому
$$
\frac{1}{\zeta(s)}\sum_{k=0}^\infty\frac{m(k)}{k^s}=\frac{n^2}{abc}-\frac{1}{n^s}  
$$
$$
+\frac{1}{a(a,n)^{s-1}}+\frac{1}{b(b,n)^{s-1}}+\frac{1}{c(c,n)^{s-1}}  
$$
$$
-\frac{n}{bc(b,c,n)^{s-1}}-\frac{n}{ac(a,c,n)^{s-1}}-\frac{n}{ab(a,b,n)^{s-1}}.   
$$
и из формул (9) и (10) вытекают аналогичные соотношения для степенных сумм $p(k)$
$$
\frac{\sum_{k=0}^{n-1}p(k)q^k}{q^n-1}=\frac{n^3}{abc}\frac{1}{q^n-1}-\frac{1}{q-1}      
$$
$$
+\frac{n}{a}\frac{1}{q^{\frac{n}{(a,n)}}-1}+\frac{n}{b}\frac{1}{q^{\frac{n}{(b,n)}}-1}+\frac{n}{c}\frac{1}{q^{\frac{n}{(c,n)}}-1}   
$$
$$
-\frac{n^2}{bc}\frac{1}{q^{\frac{n}{(b,c,n)}}-1}-\frac{n^2}{ac}\frac{1}{q^{\frac{n}{(a,c,n)}}-1}-\frac{n^2}{ab}\frac{1}{q^{\frac{n}{(a,b,n)}}-1},   
$$
$$
\frac{1}{\zeta(s)}\sum_{k=0}^\infty\frac{p(k)}{k^s}=\frac{1}{abcn^{s-3}}-1
$$
$$
+\frac{(a,n)^s}{an^{s-1}}+\frac{(b,n)^s}{bn^{s-1}}+\frac{(c,n)^s}{cn^{s-1}}  
$$
$$
-\frac{(b,c,n)^s}{bcn^{s-2}}-\frac{(a,c,n)^s}{acn^{s-2}}-\frac{(a,b,n)^s}{abn^{s-2}}.   
$$

Согласно \cite{7} характеристический многочлен монодромии квазиоднородной особенности функции трёх переменных с регулярной системой весов $(a,b,c;n)$ равен
$$
\frac{(1-q^n)^{2g-2+r}\prod_{d|n,d\in\{a,b,c\}}(1-q^{n/d})}{(1-q)\prod_{\alpha_i|n}(1-q^{n/\alpha_i})},   
$$
где $\{g,(\alpha_1,\beta_1),\dots,(\alpha_r,\beta_r)\}$ -- инварианты Зейферта особенности, определяемые $\mathbb C^*$-действием. Используя формулы (14) и (15) для $a(k)=m(k)$ и $a(k)=p(k)$, получим
$$
\frac{\sum_{k=0}^{n-1}m(k)q^k}{q^n-1}=\frac{2g-2+r}{q-1}-\frac{1}{q^n-1}
+\sum\limits_{\substack{d|n\\d\in\{a,b,c\}}}\frac{1}{q^d-1}-\sum_{\alpha_i|n}\frac{1}{q^{\alpha_i}-1},   
$$
$$
\frac{1}{\zeta(s)}\sum_{k=0}^\infty\frac{m(k)}{k^s}=2g-2+r-\frac{1}{n^s}
+\sum\limits_{\substack{d|n\\d\in\{a,b,c\}}}\frac{1}{d^s}-\sum_{\alpha_i|n}\frac{1}{\alpha_i^s},
$$
$$
\frac{\sum_{k=0}^{n-1}p(k)q^k}{q^n-1}=\frac{n(2g-2+r)}{q^n-1}-\frac{1}{q-1}
+\sum\limits_{\substack{d|n\\d\in\{a,b,c\}}}\frac{n/d}{q^{n/d}-1}-\sum_{\alpha_i|n}\frac{n/\alpha_i}{q^{n/\alpha_i}-1},   
$$
$$
\frac{n^{s-1}}{\zeta(s)}\sum_{k=0}^\infty\frac{p(k)}{k^s}=2g-2+r-n^{s-1}+\sum\limits_{\substack{d|n\\d\in\{a,b,c\}}}d^{s-1}-\sum_{\alpha_i|n}\alpha_i^{s-1}.
$$

\bigskip

\section{Простые и параболические особенности}

\bigskip

  В этом и следующем параграфах для перечисленных особенностей мы выпишем представления из формулы (15) для рациональных функций
$$
\frac{\sum_{k=0}^{n-1}m(k)q^k}{1-q^n}, \frac{\sum_{k=0}^{n-1}p(k)q^k}{1-q^n}.
$$
Для простых особенностей $A_l$, $D_l$, $E_6$, $E_7$, $E_8$
$$
\sum_{k=0}^{n-1}m(k)q^k=\sum_j q^{m_j}, 
$$
числа $m_j$ совпадают с показателями Кокстера соответствующих конечных групп отражений.

$A_l$, $n=l+1$,
$$
\frac{1}{1-q}-\frac{1}{1-q^n}, 
$$
$$
-\frac{1}{1-q}+\frac{n}{1-q^n}. 
$$
$D_l$, $n=2l-2$,
$$
\frac{1}{1-q}-\frac{1}{1-q^2}+\frac{1}{1-q^{n/2}}-\frac{1}{1-q^n}, 
$$
$$
-\frac{1}{1-q}+\frac{2}{1-q^2}-\frac{n/2}{1-q^{n/2}}+\frac{n}{1-q^n}. 
$$
$E_6$, $n=12$,
$$
\frac{1}{1-q}-\frac{1}{1-q^2}-\frac{1}{1-q^3}+\frac{1}{1-q^4}+\frac{1}{1-q^6}-\frac{1}{1-q^{12}}, 
$$
$$
-\frac{1}{1-q}+\frac{2}{1-q^2}+\frac{3}{1-q^3}-\frac{4}{1-q^4}-\frac{6}{1-q^6}+\frac{12}{1-q^{12}}. 
$$
$E_7$, $n=18$,
$$
\frac{1}{1-q}-\frac{1}{1-q^2}-\frac{1}{1-q^3}+\frac{1}{1-q^6}+\frac{1}{1-q^9}-\frac{1}{1-q^{18}}, 
$$
$$
-\frac{1}{1-q}+\frac{2}{1-q^2}+\frac{3}{1-q^3}-\frac{6}{1-q^6}-\frac{9}{1-q^9}+\frac{18}{1-q^{18}}. 
$$
$E_8$, $n=30$,
$$
\frac{1}{1-q}-\frac{1}{1-q^2}-\frac{1}{1-q^3}-\frac{1}{1-q^5}+\frac{1}{1-q^6}+\frac{1}{1-q^{10}}+\frac{1}{1-q^{15}}-\frac{1}{1-q^{30}}, 
$$
$$
-\frac{1}{1-q}+\frac{2}{1-q^2}+\frac{3}{1-q^3}+\frac{5}{1-q^5} 
-\frac{6}{1-q^6}-\frac{10}{1-q^{10}}-\frac{15}{1-q^{15}}+\frac{30}{1-q^{30}}. 
$$
$\tilde E_6=P_8$, $n=3$,
$$
\frac{3}{1-q}-\frac{1}{1-q^3},
$$
$$
-\frac{1}{1-q}+\frac{9}{1-q^3}. 
$$
$\tilde E_7=X_9$, $n=4$,
$$
\frac{2}{1-q}+\frac{1}{1-q^2}-\frac{1}{1-q^4},
$$
$$
-\frac{1}{1-q}-\frac{2}{1-q^2}+\frac{8}{1-q^4}. 
$$
$\tilde E_8=J_{10}$, $n=6$,
$$
\frac{1}{1-q}+\frac{1}{1-q^2}+\frac{1}{1-q^3}-\frac{1}{1-q^6},
$$
$$
-\frac{1}{1-q}+\frac{2}{1-q^2}+\frac{3}{1-q^4}+\frac{6}{1-q^6}. 
$$

\bigskip

\section{Исключительные особенности Арнольда}

\bigskip

$U_{12}$, $n=12$,
$$
\frac{1}{1-q}+\frac{1}{1-q^3}-\frac{1}{1-q^4}-\frac{1}{1-q^{12}}, 
$$
$$
-\frac{1}{1-q}-\frac{3}{1-q^3}+\frac{4}{1-q^4}+\frac{12}{1-q^{12}}. 
$$
$S_{12}$, $n=13$,
$$
\frac{1}{1-q}-\frac{1}{1-q^{13}}, 
$$
$$
-\frac{1}{1-q}+\frac{n}{1-q^{13}}. 
$$
$S_{11}$, $n=16$,
$$
\frac{1}{1-q}-\frac{1}{1-q^2}+\frac{1}{1-q^4}-\frac{1}{1-q^{16}}, 
$$
$$
-\frac{1}{1-q}+\frac{4}{1-q^4}-\frac{8}{1-q^8}+\frac{16}{1-q^{16}}. 
$$
$Q_{12}$, $n=15$,
$$
\frac{1}{1-q}-\frac{1}{1-q^3}+\frac{1}{1-q^5}-\frac{1}{1-q^{15}}, 
$$
$$
-\frac{1}{1-q}+\frac{3}{1-q^3}-\frac{5}{1-q^5}+\frac{15}{1-q^{15}}. 
$$
$Q_{11}$, $n=18$,
$$
\frac{1}{1-q}-\frac{1}{1-q^2}+\frac{1}{1-q^6}-\frac{1}{1-q^{18}}, 
$$
$$
-\frac{1}{1-q}+\frac{3}{1-q^3}-\frac{9}{1-q^9}+\frac{18}{1-q^{18}}. 
$$
$Q_{10}$, $n=24$,
$$
\frac{1}{1-q}-\frac{1}{1-q^2}-\frac{1}{1-q^3}+\frac{1}{1-q^6}+\frac{1}{1-q^8}-\frac{1}{1-q^{24}}, 
$$
$$
-\frac{1}{1-q}+\frac{3}{1-q^3}+\frac{4}{1-q^4}-\frac{8}{1-q^8}-\frac{12}{1-q^{12}}+\frac{24}{1-q^{24}}. 
$$
$W_{13}$, $n=16$,
$$
\frac{1}{1-q}-\frac{1}{1-q^4}+\frac{1}{1-q^8}-\frac{1}{1-q^{16}}, 
$$
$$
-\frac{1}{1-q}+\frac{2}{1-q^2}-\frac{4}{1-q^4}+\frac{16}{1-q^{16}}. 
$$
$W_{12}$, $n=20$,
$$
\frac{1}{1-q}-\frac{1}{1-q^2}+\frac{1}{1-q^4}-\frac{1}{1-q^5}+\frac{1}{1-q^{10}}-\frac{1}{1-q^{20}}, 
$$
$$
-\frac{1}{1-q}+\frac{2}{1-q^2}-\frac{4}{1-q^4}+\frac{5}{1-q^5}-\frac{10}{1-q^{10}}+\frac{20}{1-q^{20}}. 
$$
$Z_{13}$, $n=18$,
$$
\frac{1}{1-q}-\frac{1}{1-q^3}+\frac{1}{1-q^9}-\frac{1}{1-q^{18}}, 
$$
$$
-\frac{1}{1-q}+\frac{2}{1-q^2}-\frac{6}{1-q^6}+\frac{18}{1-q^{18}}. 
$$
$Z_{12}$, $n=22$,
$$
\frac{1}{1-q}-\frac{1}{1-q^2}+\frac{1}{1-q^{11}}-\frac{1}{1-q^{22}}, 
$$
$$
-\frac{1}{1-q}+\frac{2}{1-q^2}-\frac{11}{1-q^{11}}+\frac{22}{1-q^{22}}. 
$$
$Z_{11}$, $n=30$,
$$
\frac{1}{1-q}-\frac{1}{1-q^2}-\frac{1}{1-q^3}+\frac{1}{1-q^6}+\frac{1}{1-q^{15}}-\frac{1}{1-q^{30}}, 
$$
$$
-\frac{1}{1-q}+\frac{2}{1-q^2}+\frac{5}{1-q^5}-\frac{10}{1-q^{10}}-\frac{15}{1-q^{15}}+\frac{30}{1-q^{30}}. 
$$
$E_{14}$, $n=24$,
$$
\frac{1}{1-q}-\frac{1}{1-q^3}-\frac{1}{1-q^4}+\frac{1}{1-q^8}+\frac{1}{1-q^{12}}-\frac{1}{1-q^{24}}, 
$$
$$
-\frac{1}{1-q}+\frac{2}{1-q^2}+\frac{3}{1-q^3}-\frac{6}{1-q^6}-\frac{8}{1-q^8}+\frac{24}{1-q^{24}}. 
$$
$E_{13}$, $n=30$,
$$
\frac{1}{1-q}-\frac{1}{1-q^2}-\frac{1}{1-q^5}+\frac{1}{1-q^{10}}+\frac{1}{1-q^{15}}-\frac{1}{1-q^{30}}, 
$$
$$
-\frac{1}{1-q}+\frac{2}{1-q^2}+\frac{3}{1-q^3}-\frac{6}{1-q^6}-\frac{15}{1-q^{15}}+\frac{30}{1-q^{30}}. 
$$
$E_{12}$, $n=42$,
$$
\frac{1}{1-q}-\frac{1}{1-q^2}-\frac{1}{1-q^3}+\frac{1}{1-q^6}
-\frac{1}{1-q^7}+\frac{1}{1-q^{14}}+\frac{1}{1-q^{21}}-\frac{1}{1-q^{42}}, 
$$
$$
-\frac{1}{1-q}+\frac{2}{1-q^2}+\frac{3}{1-q^3}-\frac{6}{1-q^6}
+\frac{7}{1-q^7}-\frac{14}{1-q^{14}}-\frac{21}{1-q^{21}}+\frac{42}{1-q^{42}}. 
$$

\end {document}